\documentclass{amsart}
\usepackage{amssymb, mathrsfs, amsmath, verbatim}

\theoremstyle{plain}
\newtheorem{thm}{Theorem}[section]
\newtheorem{lemma}[thm]{Lemma}
\newtheorem{prop}[thm]{Proposition}
\newtheorem{cor}[thm]{Corollary}

\newtheorem*{dimineq}{The Dimension Inequality}
\newtheorem*{nonneg}{Non-negativity}
\newtheorem*{vanishing}{Vanishing}
\newtheorem*{positivity}{Positivity}
\newtheorem*{gradeconj}{The Grade Conjecture}
\newtheorem*{intthm}{The Intersection Theorem}

\theoremstyle{remark}
\newtheorem{remark}[thm]{Remark}

\newtheoremstyle{TheoremNum}
    {\topsep}{\topsep}              
    {\itshape}                      
    {}                              
    {\bfseries}                     
    {.}                             
    { }                             
    {\thmname{#1}\thmnote{ \bfseries #3}}
\theoremstyle{TheoremNum}
\newtheorem{thmn}{Theorem}

\DeclareMathOperator{\pd}{pd}
\DeclareMathOperator{\grade}{grade}
\DeclareMathOperator{\height}{ht}
\DeclareMathOperator{\depth}{depth}
\DeclareMathOperator{\codim}{codim}

\DeclareMathOperator{\ann}{ann}
\DeclareMathOperator{\Tor}{Tor}
\DeclareMathOperator{\Hom}{Hom}
\DeclareMathOperator{\Ext}{Ext}
\DeclareMathOperator{\Ass}{Ass}
\DeclareMathOperator{\coker}{coker}
\DeclareMathOperator{\Supp}{Supp}

\DeclareMathOperator{\Tot}{Tot}
\DeclareMathOperator{\rank}{rank}

\newcommand{\frob}[1]{{^{#1}\!\!}}
\renewcommand{\i}[1]{\mathfrak{#1}}
\newcommand{\p}{\i{p}}
\newcommand{\q}{\i{q}}
\newcommand{\m}{\i{m}}

\newcommand{\ZZ}{\mathbb{Z}}
\newcommand{\lto}{\mathop{\longrightarrow\,}\limits}

\title{The Grade Conjecture and Asymptotic Intersection Multiplicity}
\author{Jesse Beder}

\begin{document}

\maketitle

\begin{abstract}
Given a finitely generated module $M$ over a local ring $A$ of characteristic $p$ with $\pd M < \infty$, we study the asymptotic intersection multiplicity $\chi_\infty(M, A/\underline{x})$, where $\underline{x} = (x_1, \ldots, x_r)$ is a system of parameters for $M$. We show that there exists a system of parameters such that $\chi_\infty$ is positive if and only if $\dim \Ext^{d-r}(M, A) = r$, where $d = \dim A$ and $r = \dim M$. We use this to prove several results relating to the Grade Conjecture, which states that $\grade M + \dim M = \dim A$ for any module $M$ with $\pd M < \infty$.
\end{abstract}

\section{Introduction}

In 1965, Serre \cite{serre_localalg} proved the following theorem about intersection multiplicities:

\begin{thm}[Serre]
Let $A$ be a regular local ring either containing a field or unramified over a discrete valuation ring, and let $M$ and $N$ be finitely generated $A$-modules with $l(M\otimes N) < \infty$. Then, writing
\[ \chi(M, N) = \sum_{i=0}^{\dim A} (-1)^i l(\Tor_i(M, N)), \]
we have
\begin{enumerate}
\item $\dim M + \dim N \leq \dim A$
\item $\chi(M, N) \geq 0$
\item $\chi(M, N) = 0$ if and only if $\dim M + \dim N < \dim A$.
\end{enumerate}
\end{thm}

He also proved (1) for an arbitrary regular local ring, and conjectured that (2) and (3) hold as well. These statements can be further generalized to an arbitrary local ring when one of the modules has finite projective dimension, and they are known as:

\begin{dimineq} $\dim M + \dim N \leq \dim A$; \end{dimineq}
\begin{nonneg} $\chi(M, N) \geq 0$; \end{nonneg}
\begin{vanishing} $\chi(M, N) = 0$ if $\dim M + \dim N < \dim A$; and \end{vanishing}
\begin{positivity} $\chi(M, N) > 0$ if $\dim M + \dim N = \dim A$. \end{positivity}

When $A$ is regular local, nonnegativity was proved by Gabber \cite{gabber_nonneg}, and vanishing was proved by Roberts \cite{roberts_vanishing} and Gillet-Soul\'{e} \cite{gs_vanishing} independently; positivity is still open when $A$ is a ramified regular local ring. When $A$ is not regular, these conjectures are mostly open. Nonnegativity is false in general \cite{dutta_etal}.

Positivity is even unknown in the following special case: if $\pd M < \infty$ and $x_1, \ldots, x_r$ is a system of parameters for $M$, is $\chi(M, A/\underline{x}) > 0$?

In 1982, Dutta \cite{dutta_frobmult} introduced an asymptotic multiplicity $\chi_\infty$, defined in Section~\ref{notation}, to investigate vanishing and positivity over a local ring of characteristic $p$. In particular, he showed that $\chi_\infty(M, N) = 0$ if $\dim M + \dim N < \dim A$, and that $\chi_\infty(M, N) > 0$ if $\dim M + \dim N = \dim A$ and $M$ is Cohen-Macaulay.

However, given a module $M$ of finite projective dimension, and $x_1, \ldots, x_r$ a system of parameters for $M$ \emph{that is also an $A$-sequence}, it is unknown whether $\chi_\infty(M, A/\underline{x}) > 0$; in this case, we do know that $\chi(M, A/\underline{x}) > 0$, due to Serre \cite{serre_localalg} and Lichtenbaum \cite{lichtenbaum}.

We will consider $\chi_\infty(M, A/\underline{x})$ when $x_1, \ldots, x_r$ is a system of parameters for $M$. We prove the following:

\begin{thmn}[\ref{chi_ext_dim}]
Suppose $\pd M < \infty$, where $d = \dim A$ and $r = \dim M$. Then there is a system of parameters $x_1, \ldots, x_r$ for $M$, that is part of a system of parameters for $A$, such that
\[ \chi_\infty(M, A/\underline{x}) > 0 \mbox{ if and only if } \dim \Ext^{d-r}(M, A) = r. \]
\end{thmn}

We will also show a special case of asymptotic positivity:

\begin{thmn}[\ref{chi_pos_lowpd}]
Let $d = \dim A$ and $r = \dim M$, and suppose that $\pd M = d - r$. Then any system of parameters $x_1, \ldots, x_r$ for $M$ is part of a system of parameters for $A$, and
\[ \chi_\infty(M, A/\underline{x}) > 0. \]
\end{thmn}

The equivalence in Theorem \ref{chi_ext_dim} suggests that we investigate the condition
\[ \dim \Ext^{d-r}(M, A) = r, \]
and in Section \ref{grade_sec} this will lead us to study the following conjecture of Peskine and Szpiro \cite[Conjecture (f) of Chapter II]{ps_ihes}:

\begin{gradeconj}
Suppose that $\pd M < \infty$. Then
\[ \grade M + \dim M = \dim A. \]
\end{gradeconj}

The Grade Conjecture is known in some specific cases: when $M$ is perfect or $A$ is Cohen-Macaulay (these are due to Peskine and Szpiro \cite{ps_ihes}); and in the graded case, when $M = \bigoplus M_i$ is a graded module over a graded ring $A = \bigoplus A_i$ with $A_0$ artinian (this is also due to Peskine and Szpiro \cite{ps_syz}). Foxby \cite{foxby_flatcomplex} showed that the Grade Conjecture holds if $A$ is complete and equidimensional in the equicharacteristic case, and Roberts \cite{roberts_equidim} showed this is true in arbitrary characteristic.

In Section \ref{grade_sec}, we will work toward a connection between the Grade Conjecture and $\chi_\infty$, and along the way we will prove Roberts' Theorem. We will then prove the following:

\begin{thmn}[\ref{chi_grade_conj}]
Let $A$ be a local ring in characteristic $p$, and suppose that $\pd M < \infty$. Assume the Grade Conjecture holds for $M$. Then there is a system of parameters $x_1, \ldots, x_r$ for $M$, that is part of a system of parameters for $A$, such that
\[ \chi_\infty(M, A/\underline{x}) > 0. \]
\end{thmn}

This connection gives us a few new cases in which $\dim \Ext^{d-r}(M, A) = r$ in characteristic $p$ (recall $d = \dim A$ and $r = \dim M$): when $\pd M = d - r$, and when $\dim M = 1$ (Theorem \ref{ext_dim_lowpd} and Corollary \ref{dim_one_dim}).

Finally, in Section \ref{equichar_zero}, we use theorems of Hochster and Kurano to translate these characteristic $p$ results to equicharacteristic zero.

\subsection{Notation}\label{notation}

Unless otherwise specified, all modules are finitely generated. All $\Tor$ and $\Ext$ modules are computed over $A$ unless otherwise specified.

If $M$ is an $A$-module and $x_1, \ldots, x_r$ elements of $A$, then write $e(\underline{x}; M)$ for the multiplicity of $M$ with respect to $x_1, \ldots, x_r$. Depending on the context, we write $\underline{x}$ for either the sequence $x_1, \ldots, x_r$ or the ideal generated by that sequence.

When $A$ is of prime characteristic $p$, we let $f: A\to A$ be the Frobenius map, $x \mapsto x^p$, and let $f^n$ be its $n$th iteration. We write $\frob{f^n}A$ for the bialgebra $A$ with the action on the left by $f^n$ and the right by the identity. If $M$ is an $A$-module, then we write $F^n(M) = M\otimes\frob{f^n}A$ \cite{ps_ihes}.

When the limits exist, we define \cite{dutta_frobmult}
\[ \chi_\infty(M, N) = \lim_{n\to\infty} \frac{\chi(F^n(M), N)}{p^{n\codim M}} \]
and
\[ e_\infty(\underline{x}; M) = \lim_{n\to\infty} \frac{e(\underline{x}; F^n(M))}{p^{n\codim M}}, \]
where $\codim M = \dim A - \dim M$.

If $M$ is an $A$-module, then
\[ \grade M = \depth \ann M \]
is the length of the longest $A$-sequence contained in $\ann M$. We say $M$ is \emph{perfect} if $\pd M = \grade M$.

\section{Ext and Asymptotic Limits}\label{ext_sec}

Throughout this section, we suppose $A$ is complete and has characteristic $p$, and let $d = \dim A$.

We first need some facts about asymptotic limits \cite[Proposition 1]{seibert_frob}.

\begin{lemma}[Seibert]\label{tor_length}
Let $L_\bullet$ be a complex of finitely generated free modules with homologies of finite length. For each $i\geq 0$, the limit
\[ \lim_{n\to\infty} \frac{l(H_i(F^n(L_\bullet)))}{p^{nd}} \]
exists, and if $H_0(L_\bullet)\neq 0$, then
\[ \lim_{n\to\infty} \frac{l(H_0(F^n(L_\bullet)))}{p^{nd}} \geq 1. \]
\end{lemma}

We now prove an equivalence we will use several times.

\begin{prop}\label{ext_dim_prime}
Suppose $\pd M < \infty$, and set $r = \dim M$. Then
\[ \dim \Ext^{d-r}(M, A) = r \]
if and only if there is a prime $\p\in \Ass M$ with $\dim A/\p = r$ and $\height \p = d - r$.
\end{prop}
\begin{proof}
If $\dim\Ext^{d-r}(M, A) = r$, then there is some prime
\[ \p\in \Supp \Ext^{d-r}(M, A)\subseteq \Supp M \]
with $\dim A/\p = r$. Such a prime is necessarily minimal over $\ann M$, so it is in $\Ass M$. Furthermore, $\pd M_\p \geq d - r$ (since $\Ext^{d-r}_{A_\p}(M_\p, A_\p)\neq 0$), so
\[ \height \p = \dim A_\p \geq \depth A_\p \geq \pd M_\p \geq d-r \]
by Auslander-Buchsbaum. Since $\height \p$ certainly cannot be larger, $\height\p = d - r$.

Conversely, suppose that $\p\in \Ass M$ with $\dim A/\p = r$ and $\height \p = d - r$. Then $M_\p$ is of finite length and finite projective dimension, so by Lemma \ref{dim_zero_perfect}, $M_\p$ is perfect of projective dimension $d-r$ over $A_\p$, and hence
\[ \Ext^{d-r}(M, A)_\p \neq 0, \]
so $\dim \Ext^{d-r}(M, A) = r$ (the dimension clearly cannot be larger than $r$).
\end{proof}

Using Proposition \ref{ext_dim_prime}, we can relate $e_\infty$ to the dimension of $\Ext^{d-r}(M, A)$.

\begin{thm}\label{ext_dim}
Suppose $\pd M < \infty$, and set $r = \dim M$. Then
\[ \dim \Ext^{d-r}(M, A) = r \]
if and only if $e_\infty(\underline{x}; M) > 0$ for some (= every) system of parameters $x_1, \ldots, x_r$ for $M$ that is part of a system of parameters for $A$.
\end{thm}
\begin{proof}
We have \cite[Theorem 14.7]{mats}
\[ e(\underline{x}; F^n(M)) = \sum_{\p} e(\underline{x}; A/\p) l(F^n(M)_\p), \]
where the sum is taken over all primes $\p\in \Supp M$ with $\dim A/\p = r$. Now, due to the associativity of the tensor product, it follows \cite{ps_ihes} that
\[ l(F^n_A(M)_\p) = l(F^n_{A_\p}(M_\p)). \]
By Proposition \ref{ext_dim_prime}, $\dim \Ext^{d-r}(M, A) = r$ if and only if there is a prime $\mathfrak{p}\in \Supp M$ with $\dim A/\p = r$ and $\height \p = d-r$; which, using Lemma \ref{tor_length} on each $l(F^n_{A_\p}(M_\p))$ in the above formula, happens if and only if $e_\infty(\underline{x}; M) > 0$.
\end{proof}

\begin{remark}
If the Dimension Inequality is true, then it implies that any system of parameters for $M$ is part of a system of parameters for $A$, so the last assumption in the above theorem would be superfluous.
\end{remark}

We will now relate $e_\infty$ to $\chi_\infty$ for a suitably chosen system of parameters:

\begin{prop}\label{mult_eq_chi}
Suppose $\pd M < \infty$ and $\dim M < \dim A$. Then there is a system of parameters $x_1, \ldots, x_r$ for $M$, that is part of a system of parameters for $A$, such that
\[ e(\underline{x}; M) = \chi(M, A/\underline{x}). \]
\end{prop}
\begin{proof}
We first claim that we can choose a system of parameters $x_1, \ldots, x_r$ for $M$ that is part of a system of parameters for $A$, such that the higher Koszul homologies have finite length, i.e.,
\[ l(H_i(\underline{x}; A)) <\infty \mbox{ for all $i\geq 1$.} \]
To do so, choose $x_{i+1}\in \m$
\begin{enumerate}
\item a parameter on $M/(x_1, \ldots, x_i)M$
\item a parameter on $A/(x_1, \ldots, x_i)A$
\item not in any prime $\p\not= \m$ in $\Ass(A/(x_1, \ldots, x_i)A)$
\end{enumerate}
Then, if we localize the Koszul complex $K_\bullet(\underline{x}; A)$ at any nonmaximal prime, all the $x_i$ will be either units or nonzerodivisors mod $x_1, \ldots, x_{i-1}$, so the higher homologies will vanish, and hence they all have finite length.

Now, let $D_{\bullet\bullet}$ be the double complex from tensoring $K_\bullet(\underline{x}; A)$ with a free resolution $L_\bullet$ for $M$. There are two associated spectral sequences, both converging to $H^*(\Tot(D_{\bullet\bullet}))$:
\begin{enumerate}
\item $_I E_2^{s, t} = H_t(\underline{x}; M)$ when $s = 0$, and zero otherwise
\item $_{II} E_2^{s, t} = \Tor_s(M, H_t(\underline{x}; A))$
\end{enumerate}
Since all the above $E_2$ terms have finite length, and there are only finitely many of them (they vanish for $s > \pd M$ and $t > r$), the Euler characteristics of the two $E_2$ terms are equal, so we have
\[ \sum_t (-1)^t l(H_t(\underline{x}; M)) = \sum_{s, t} (-1)^{s+t} l(\Tor_s(M, H_t(\underline{x}; A))). \]
The left-hand side is just $e(\underline{x}; M)$ (see, e.g., \cite[Theorem 1 of Chapter IV]{serre_localalg}), and the right side is just
\[ \sum_t (-1)^t \chi(M, H_t(\underline{x}; A)). \]
For $t\geq 1$, the $H_t(\underline{x}; A)$ have finite length, by the choice of $\underline{x}$ above; since $\dim M < \dim A$, we have $\chi(M, H_t(\underline{x}; A)) = 0$ for $t \geq 1$. Thus
\[ e(\underline{x}; M) = \chi(M, H_0(\underline{x}; A)) = \chi(M, A/\underline{x}), \]
as desired.
\end{proof}

Note that the choice of $x_1, \ldots, x_r$ depends only on $\Supp M$, so the same system of parameters can be used for $F^n(M)$ for all $n$. We can therefore use Proposition \ref{mult_eq_chi} to translate Theorem \ref{ext_dim} to $\chi_\infty$:

\begin{thm}\label{chi_ext_dim}
Suppose $\pd M < \infty$, and set $d = \dim A$ and $r = \dim M$. Then there is a system of parameters $x_1, \ldots, x_r$ for $M$, that is part of a system of parameters for $A$, such that
\[ \chi_\infty(M, A/\underline{x}) > 0 \mbox{ if and only if } \dim \Ext^{d-r}(M, A) = r. \]
\end{thm}

For the proof of the next theorem, we need the following two results. First, the following theorem of Dutta's \cite{dutta_extfrobII}:
\begin{thm}[Dutta]\label{dutta_hom}
Let $F_\bullet$ be a complex of finitely-generated free modules with homologies of finite length. Let $N$ be a finitely generated module. Let $W_{jn}$ denote the $j$th homology of $\Hom(F^n(F_\bullet), N)$, and write
\[ \tilde{N} = \Hom(H_\m^d(N), E), \]
where $E = E(k)$ is the injective hull of $k$.

We have the following:
\begin{enumerate}
\item If $\dim N < \dim A$, then $\lim l(W_{jn})/p^{nd} = 0$.
\item If $\dim N = \dim A$, and
\begin{enumerate}
\item $j < d$, then $\lim l(W_{jn})/p^{nd} = 0$;
\item $j = d$, then $\lim l(W_{jn})/p^{nd} = \lim l(F^n(H_0(F_\bullet))\otimes \tilde{N})/p^{nd}$, which is positive;
\item $j > d$, then $\lim l(W_{jn})/p^{nd} = \lim l(H_{j-d}(F^n(F_\bullet))\otimes \tilde{N})/p^{nd}$.
\end{enumerate}
\end{enumerate}
\end{thm}

We will also use Peskine and Szpiro's Intersection Theorem:

\begin{intthm}
Suppose that $\pd M < \infty$, and let $N$ be another finitely generated $A$-module such that $l(M\otimes N) < \infty$. Then
\[ \dim N \leq \pd M. \]
\end{intthm}
Peskine and Szpiro \cite{ps_ihes} first proved this when $A$ has characteristic $p$, and in many equicharacteristic cases; Hochster \cite{hochster_equichar} proved the equicharacteristic case; and Roberts \cite{roberts_intthm} proved the mixed characteristic case.

We will now show a special case of asymptotic positivity:

\begin{thm}\label{chi_pos_lowpd}
Let $d = \dim A$ and $r = \dim M$, and suppose that $\pd M = d - r$. Then any system of parameters $x_1, \ldots, x_r$ for $M$ is part of a system of parameters for $A$, and
\[ \chi_\infty(M, A/\underline{x}) > 0. \]
\end{thm}

\begin{proof}
We write $B = A/\underline{x}$. By the Intersection Theorem,
\[ \dim B \leq \pd M = d - r, \]
and so we can choose $y_1, \ldots, y_{d-r}$ such that $l(B/\underline{y}B) < \infty$. Since $d = \dim A$, it follows that $x_1, \ldots, x_r, y_1, \ldots, y_{d-r}$ is a system of parameters for $A$, proving the first claim.

Now we let $L_\bullet$ be a free resolution for $M$ over $A$, and let $F_\bullet = \Hom(L_\bullet, B)$. We will apply Theorem \ref{dutta_hom} to $F_\bullet$ over the ring $B$, with $N = B$. We note that
\[ W_{jn} = \Tor_{d-r-j}(F^n(M), B), \]
since the length of the complex $L_\bullet$ is $\pd M = d - r$ (we are flipping the indices for $F_\bullet$; that is, $F_i = \Hom(L_{d-r-i}, B)$). By Theorem \ref{dutta_hom},
\[ \lim_{n\to\infty} \frac{l(W_{jn})}{p^{n(d-r)}} = 0, \]
for all $j < d - r$, which means that
\[ \chi_\infty(M, B) = \lim_{n\to\infty} \frac{l(F^n(M)\otimes B)}{p^{n(d-r)}} > 0 \]
by Lemma \ref{tor_length}, proving the theorem.
\end{proof}

\section{The Grade Conjecture}\label{grade_sec}

In this section, $A$ is an arbitrary local ring, unless otherwise specified.

We begin with a result of Peskine and Szpiro's \cite[Lemma 4.8]{ps_ihes}:

\begin{lemma}[Peskine-Szpiro]\label{grade_bounds}
Let $M$ be a finite $A$-module. Then
\[ \depth A \leq \grade M + \dim M \leq \dim A. \]
\end{lemma}

This immediately implies the grade conjecture when $A$ is Cohen-Macaulay:

\begin{cor}\label{cm_grade_conj}
Suppose $A$ is Cohen-Macaulay. Then
\[ \grade M + \dim M = \dim A \]
for all finitely generated modules $M$.
\end{cor}

When $M$ is perfect, the grade conjecture follows immediately from the Intersection Theorem:

\begin{cor}[Peskine-Szpiro]\label{perfect}
Suppose $\pd M < \infty$ and $M$ is perfect. Then
\[ \grade M + \dim M = \dim A. \]
\end{cor}
\begin{proof}
Let $x_1, \ldots, x_r$ be a system of parameters for $M$. By the Intersection Theorem (and since $M$ is perfect),
\[ \dim A/\underline{x} \leq \pd M = \grade M. \]
Using $\dim A/\underline{x} \geq \dim A - \dim M$ gives $\dim A \leq \grade M + \dim M$; and the reverse inequality follows from Lemma \ref{grade_bounds}.
\end{proof}

Next, we have some basic results about the annihilator of modules with finite projective dimension.

\begin{lemma}\label{dim_zero_perfect}
Suppose $\pd M < \infty$ and $\dim M = 0$. Then $M$ is perfect and $A$ is Cohen-Macaulay.
\end{lemma}
\begin{proof}
Since $\dim M = 0$, $\ann M$ is $\m$-primary, so $\grade M = \depth A$. By Auslander-Buchsbaum,
\[ \pd M + \depth M = \depth A, \]
and so $\depth M = 0$ implies $\pd M = \grade M$, i.e., $M$ is perfect; and by Corollary \ref{perfect}, $A$ is Cohen-Macaulay.
\end{proof}

\begin{lemma}\label{ht_equals_grade}
Suppose $\pd M < \infty$. Then $\height \ann M = \grade M$.
\end{lemma}
\begin{proof}
Let $x_1, \ldots, x_g$ be a maximal $A$-sequence contained in $I = \ann M$ (so $g = \grade M$). Then
\[ \Hom(M, A/\underline{x}) = \Ext^g (M, A) \neq 0. \]
Choose a minimal prime $\p$ in $\Supp M\cap \Ass A/\underline{x}$. Then
\[ \Hom(M_\p, A_\p/\underline{x} A_\p) \neq 0, \]
so $\grade M_\p = \grade M = g$ and $\height \p \geq \height I$.

Since $\p\in \Ass A/\underline{x}$, $\depth A_\p = g$; by Auslander-Buchsbaum,
\[ \depth M_\p + \pd M_\p = g; \]
and hence $\pd M_\p = g = \grade M_\p$, i.e., $M_\p$ is perfect over $A_\p$. Corollary \ref{perfect} then implies that
\[ \grade M_\p + \dim M_\p = \dim A_\p. \]
Now we choose a prime $\q$ with $I\subseteq \q\subseteq \p$, minimal over $I$, with $\dim A_p/\q A_p = \dim M_\p$. We also note that $\dim A_\q \geq \depth A_\q \geq g$.

We now have
\begin{eqnarray*}
\dim A_\p & \geq & \dim A_\p/\q A_\p + \dim A_q \\
          & = & \dim M_\p + \dim A_\q \\
          & \geq & \dim M_\p + g \\
          & = & \dim A_\p,
\end{eqnarray*}
which implies that the above inequalities are equalities; in particular, $\dim A_\q = g$, which means that $\height I \leq \height \q = g = \grade M$. Since the reverse inequality is always true, we have equality.
\end{proof}

From this, we get a special case of the grade conjecture \cite{roberts_equidim}:

\begin{thm}[Roberts]\label{equidim_grade_conj}
Suppose that $A$ is complete and equidimensional and $\pd M < \infty$. Then
\[ \grade M + \dim M = \dim A. \]
\end{thm}
\begin{proof}
We will show that $\grade M + \dim M \geq \dim A$; equality then holds because the reverse inequality is just Lemma \ref{grade_bounds}.

We choose a prime $\p\in \Supp M$ with $\height \p = \height \ann M$, and then choose a minimal prime $\q\subseteq \p$ with $\height \p/\q = \height \p$. We then have
\begin{eqnarray*}
\dim A & = & \dim A/\q \\
       & = & \dim A/\p + \height \p/\q \\
       & = & \dim A/\p + \height \p \\
       & = & \dim A/\p + \grade M \\
       & \leq & \dim M + \grade M
\end{eqnarray*}
where the first equality follows because $A$ is equidimensional; the second because $A$ is complete, and hence catenary; and the fourth from Lemma \ref{ht_equals_grade}. Finally, the inequality holds since $\dim A/\p \leq \dim M$ always.
\end{proof}

We will now prove a connection between the Grade Conjecture and $\chi_\infty$.

\begin{cor}\label{ext_dim_grade_conj}
Suppose $\pd M < \infty$ and that $\grade M + \dim M = \dim A$ (i.e., the grade conjecture holds for $M$). Let $d = \dim A$ and $r = \dim M$. Then
\[ \dim \Ext^{d-r}(M, A) = r. \]
\end{cor}
\begin{proof}
Let $\p\in \Supp M$ with $\dim A/\p = r$. Then $l(M_\p) < \infty$, so by Lemma \ref{dim_zero_perfect}, $M_\p$ is perfect over the Cohen-Macaulay ring $A_\p$. We then have
\[ d - r = \grade M \leq \grade M_\p = \dim A_\p = \height \p, \]
so we are done by Proposition \ref{ext_dim_prime}.
\end{proof}

By Theorem \ref{chi_ext_dim}, this implies:

\begin{thm}\label{chi_grade_conj}
Let $A$ be a local ring in characteristic $p$, and suppose that $\pd M < \infty$. Assume the Grade Conjecture holds for $M$. Then there is a system of parameters $x_1, \ldots, x_r$ for $M$, that is part of a system of parameters for $A$, such that
\[ \chi_\infty(M, A/\underline{x}) > 0. \]
\end{thm}

We also can use $\chi_\infty$ to show special cases for which $\dim \Ext^{d-r}(M, A) = r$.

\begin{thm}\label{ext_dim_lowpd}
Let $A$ be a local ring in characteristic $p$, and let $d = \dim A$ and $r = \dim M$. Suppose that $\pd M = d - r$. Then $\dim\Ext^{d-r}(M, A) = r$.
\end{thm}

\begin{proof}
Choose $x_1, \ldots, x_r$ as in Proposition \ref{mult_eq_chi}. We then have
\[ e_\infty(\underline{x}; M) = \chi_\infty(M, A/\underline{x}) > 0, \]
by Theorem \ref{chi_pos_lowpd}, so the theorem follows from Theorem \ref{ext_dim}.
\end{proof}

\begin{cor}\label{dim_one_dim}
Let $A$ be a local ring in characteristic $p$ of dimension $d$, and suppose $\pd M < \infty$ and $\dim M = 1$. Then $\dim \Ext^{d-1}(M, A) = 1$.
\end{cor}
\begin{proof}
By the Intersection Theorem, either $\pd M = d - 1$ or $\pd M = d$.

In the former case, the result follows from Theorem \ref{ext_dim_lowpd}. In the latter case, $A$ is Cohen-Macaulay, so the grade conjecture holds (Corollary \ref{cm_grade_conj}) and the result follows from Corollary \ref{ext_dim_grade_conj}.
\end{proof}

\section{The Equicharacteristic Zero Case}\label{equichar_zero}

One can show Theorem \ref{ext_dim_lowpd} and Corollary \ref{dim_one_dim} in equicharacteristic zero by using standard techniques to reduce to characteristic $p$. We will sketch the proof for Theorem \ref{ext_dim_lowpd}, and the same techniques can be used for Corollary \ref{dim_one_dim}, as well as the Grade Conjecture itself.

In this section, if $f_1, \ldots, f_N\in \ZZ[X_1, \ldots, X_m]$, and $x_1, \ldots, x_m$ are elements of some ring $R$, then we say that $\underline{x}$ is a solution of $f_1, \ldots, f_N$ if
\[ f_1(\underline{x}) = 0, \ldots, f_N(\underline{x}) = 0. \]
The basis for the reduction is Hochster's ``metatheorem'':

\begin{thm}[Hochster]
Let
\[ f_1, \ldots, f_N\in \ZZ[X_1, \ldots, X_d, U_1, \ldots, U_t], \]
where $X_i$ and $U_i$ are indeterminates. If $f_1, \ldots, f_N$ have a solution $(\underline{x}, \underline{u})$ in a local ring $A$ of equicharacteristic zero, with $x_1, \ldots, x_d$ forming a system of parameters for $A$, then $f_1, \ldots, f_N$ have a solution $(\underline{x}', \underline{u}')$ in a local ring $A'$ of characteristic $p$, with $x_1', \ldots, x_d'$ forming a system of parameters for $A'$.
\end{thm}

We will use a more generalized version \cite{kurano_int}:

\begin{thm}[Kurano]\label{reduction}
Let
\[ f_1, \ldots, f_N\in \ZZ[Y_1, \ldots, Y_n, X_1, \ldots, X_d, G_1, \ldots, G_l, U_1, \ldots, U_t], \]
where $X_i$, $Y_i$, $G_i$, and $U_i$ are indeterminates. If $f_1, \ldots, f_N$ have a solution $(\underline{y}, \underline{x}, \underline{g}, \underline{u})$ in a regular local ring $R$ of equicharacteristic zero, with $y_1, \ldots, y_n$ forming a regular system of parameters for $R$, and $x_1, \ldots, x_d$ forming a system of parameters for $R/\underline{g}$, then $f_1, \ldots, f_N$ have a solution $(\underline{y}', \underline{x}', \underline{g}', \underline{u}')$ in a regular local ring $R'$ of characteristic $p$, with $y_1', \ldots, y_n'$ forming a regular system of parameters for $R'$, and $x_1', \ldots, x_d'$ forming a system of parameters for $R'/\underline{g}'$.
\end{thm}

The goal, then, is to find equations that preserve the particular properties that we're interested in. The main result we need is that we can preserve the dimension of homology modules in the following sense.

\begin{lemma}\label{twomaps}
Let $R$ be a regular local ring, and let $A = R/(g_1, \ldots, g_l)$ be a quotient. Suppose that
\[ A^a \lto^{\phi} A^b \lto^{\psi} A^c \]
is a complex of free $A$-modules with homology of dimension $h$ (where we set $h = -\infty$ if the complex is exact). We choose matrices $(\overline{r}_{ij})$ and $(\overline{s}_{ij})$ that represent the maps $\phi$ and $\psi$, with $r_{ij}, s_{ij}\in R$. Then there are polynomials $f_1, \ldots, f_N$ with coefficients in $\ZZ$ in indeterminates
\begin{enumerate}
\item $Y_1, \ldots, Y_n$
\item $G_1, \ldots, G_l$
\item $U_{ij}$ and $V_{ij}$, corresponding to the matrices $r_{ij}$ and $s_{ij}$
\item $W_1, \ldots, W_t$ (for some sufficiently large $t$)
\end{enumerate}
such that
\begin{enumerate}
\item There are $\underline{y}$ and $\underline{w}$ in $R$ such that $(\underline{y}, \underline{g}, (r_{ij}), (s_{ij}), \underline{w})$ is a solution of $f_1, \ldots, f_N$, and $y_1, \ldots, y_n$ forms a regular system of parameters for $R$.
\item If $(\underline{y}', \underline{g}', (r_{ij}'), (s_{ij}'), \underline{w}')$ is a solution of $f_1, \ldots, f_N$ in a regular local ring $R'$ with $y_1', \ldots, y_n'$ forming a regular system of parameters for $R'$, then, setting $A' = R'/\underline{g}'$, and letting $\phi'$ and $\psi'$ be maps corresponding to the matrices $(\overline{r}_{ij}')$ and $(\overline{s}_{ij}'$,
\[ A'^a \lto^{\phi'} A'^b \lto^{\psi'} A'^c \]
is a complex of free $A'$ modules with homology of dimension $h$.
\end{enumerate}
\end{lemma}
The proof is very similar to the proof of \cite[Lemma 3.11]{kurano_int}. We will give an outline.
\begin{proof}(Sketch)
Let $K_\phi, I_\phi, C_\phi; K_\psi, I_\psi, C_\psi$ be the kernel, image, and cokernel of $\phi$ and $\psi$, respectively. Let $\mathbf{K}_{\phi\bullet}, \mathbf{I}_{\phi\bullet}, \mathbf{C}_{\phi\bullet}, \mathbf{K}_{\psi\bullet},  \mathbf{I}_{\psi\bullet}, \mathbf{C}_{\psi\bullet}$ be minimal free resolutions for each over $R$. This gives us exact sequences of acyclic complexes of free $R$-modules
\[
\begin{array}{ccccccccc}
0 &\to& \mathbf{K}_{\phi\bullet} &\to& \mathbf{F}_\bullet &\to& \mathbf{I}_{\phi\bullet} &\to& 0 \\
0 &\to& \mathbf{I}_{\phi\bullet} &\to& \mathbf{F}_\bullet' &\to& \mathbf{C}_{\phi\bullet} &\to& 0 \\
0 &\to& \mathbf{K}_{\psi\bullet} &\to& \mathbf{F}_\bullet'' &\to& \mathbf{I}_{\psi\bullet} &\to& 0 \\
0 &\to& \mathbf{I}_{\psi\bullet} &\to& \mathbf{F}_\bullet''' &\to& \mathbf{C}_{\psi\bullet} &\to& 0
\end{array}
\]
where $\mathbf{F}_\bullet$ is a finite free resolution of $A^a$, $\mathbf{F}_\bullet'$ and $\mathbf{F}_\bullet''$ are of $A^b$, and $\mathbf{F}_\bullet'''$ is of $A^c$.

We can similarly write down exact sequences of complexes relating the finite free resolutions $\mathbf{F}_\bullet$, etc., with a minimal free resolution of $A$ (summed $a, b$, or $c$ times, as appropriate).

Next, we consider the exact sequence
\[ 0\to I_\phi \to K_\psi \to H \to 0 \]
where $H$ is the homology of the original sequence. We extend the map $I_\phi\to K_\psi$ to a map on the complexes $\mathbf{I}_{\phi\bullet}\to \mathbf{K}_{\psi\bullet}$. This allows us to construct a free resolution for $H$ that ends with
\begin{equation}\label{coker_homology}
({I_{\phi}})_0 \oplus ({K_{\psi}})_1 \to ({K_{\psi}})_0 \to H \to 0.
\end{equation}

Now, we would like to find equations to preserve the exactness of all the complexes above, and also the dimension of the cokernel $H$ of the map in (\ref{coker_homology}).

To preserve the exactness of a finite complex of free $R$-modules, e.g.,
\[ \mathbf{F}_\bullet: 0 \to F_h \lto^{d_h} F_{h-1} \lto^{d_{h-1}} \cdots \lto^{d_2} F_1 \lto^{d_1} F_0, \]
we first choose matrices $(c_{ij}^k)$ to represent the maps $d_k$, and choose corresponding variables $C_{ij}^k$. From each pairwise composition $d_{k+1} d_k = 0$, we get a system of equations on the $C_{ij}^k$s that preserves the fact that $\mathbf{F}_\bullet$ is a complex.

To preserve the fact that $\mathbf{F}_\bullet$ is acyclic, we use the Buchsbaum-Eisenbud criteria for exactness \cite{be_exact}:
\begin{itemize}
\item[a)] To ensure the rank of each $d_k$ is at least $r_k = \sum_{i=k}^h (-1)^i \rank F_i$, we add equations so $r_k+1$ minors of the matrices $C_{ij}^k$ vanish.
\item[b)] To ensure that the grade ($=$ height, since $R$ is regular local) of the determinental ideal $I_{r_k}(d_k)$ is at least $k$, we first add $n - \height I$ elements of $R$ to $I$ to reduce to the case where $I$ is $\m$-primary, and then write equations such that $\sqrt{I} = \m = (y_1, \ldots, y_n)$.
\end{itemize}

Finally, preserving the dimension of the cokernel of a map of free $R$-modules is similar (this is preserving the height of an ideal, which is a bit more work than above, where we prevent the height from decreasing), which allows us to preserve the dimension of the homology $H$.
\end{proof}

This allows us to preserve many properties of modules. In particular:

\begin{lemma}\label{preserve_module}
Let $R$ be a regular local ring, let $A = R/(g_1, \ldots, g_l)$ be a quotient, and let $M$ be an $A$-module of finite projective dimension. We let $d = \dim A, r = \dim M, g = \grade_A M$, and $e_i = \dim \Ext^i_A(M, A)$. Let
\[ F_\bullet:\qquad 0\lto A^{b_h}\lto^{\phi^h} A^{b_{h-1}} \lto \cdots \lto^{\phi^1} A^{b_0} \]
be a minimal free resolution for $M$. We choose matrices $(\overline{r}_{ij}^k)$ that represent the maps $\phi^k$, for $1\leq k\leq h$. Then there are polynomials $f_1, \ldots, f_N$ with coefficients in $\ZZ$ in indeterminates
\begin{enumerate}
\item $Y_1, \ldots, Y_n$
\item $X_1, \ldots, X_d$
\item $G_1, \ldots, G_l$
\item $U_{ij}^k$, corresponding to the matrices $r_{ij}^k$
\item $W_1, \ldots, W_t$ (for some sufficiently large $t$)
\end{enumerate}
such that
\begin{enumerate}
\item There are $\underline{y}, \underline{x}, \underline{w}$ in $R$ such that $(\underline{y}, \underline{x}, \underline{g}, (r_{ij}^k), \underline{w})$ is a solution to $f_1, \ldots, f_N$, with $y_1, \ldots, y_n$ a regular system of parameters for $R$ and $x_1, \ldots, x_d$ a system of parameters for $A$.
\item If $(\underline{y}', \underline{x}', \underline{g}', ({r'}_{ij}^k), \underline{w}')$ is a solution of $f_1, \ldots, f_N$ in a regular local ring $R'$ with $y_1', \ldots, y_n'$ forming a regular system of parameters for $R'$ and $x_1', \ldots, x_d'$ is a system of parameters for $A' = R'/\underline{g}'$, then letting ${\phi'}^k$ be maps corresponding to the matrices $(\overline{r'}_{ij}^k)$, we have
\[ F_\bullet':\qquad 0\lto A'^{b_h}\lto^{{\phi'}^h} A'^{b_{h-1}} \lto \cdots \lto^{{\phi'}^1} A'^{b_0} \]
is exact; and, setting $M' = \coker({\phi'}^1)$, we have $r = \dim M', g = \grade_{A'} M'$, and $e_i = \dim \Ext^i_{A'}(M', A')$.
\end{enumerate}
\end{lemma}

\begin{proof}
For each pair of maps in the free resolution for $M$,
\[ A^{b_{k+1}} \lto^{\phi^{k+1}} A^{b_k} \lto^{\phi^k} A^{b_{k-1}}, \]
we choose equations as in Lemma \ref{twomaps}, using the same variables $U_{ij}^k$ for each map $\phi^k$. This preserves the projective resolution of $M$, as well as $\dim M$.

Now, taking $\Hom(\--, A)$ of the free resolution for $M$, for each pair of maps
\[ A^{b_{k-1}} \lto^{{\phi^k}^*} A^{b_k} \lto^{{\phi^{k+1}}^*} A^{b_{k+1}}, \]
we choose equations as in Lemma \ref{twomaps}, using the same variables for the dual of a map as we used for the map itself. This preserves $\dim \Ext^k_A(M, A)$ (including whether $\Ext^k_A(M, A) = 0$), and so it also preserves $\grade M$. Combining these sets of equations, we can preserve all the numerical data, as desired.
\end{proof}

Putting this together with Theorem \ref{reduction}, we can now extend Theorem \ref{ext_dim_lowpd} to equicharacteristic:

\begin{thm}
Let $A$ be a local ring in equicharacteristic, and let $d = \dim A$ and $r = \dim M$. Suppose that $\pd M = d - r$. Then $\dim\Ext^{d-r}(M, A) = r$.
\end{thm}

\begin{proof}
Suppose there is a counterexample, that is, a local ring $A$ of equicharacteristic zero (since the theorem is true in characteristic $p$) and an $A$-module $M$ of with $\pd M = d - r$ and $\dim\Ext^{d-r}(M, A) \neq r$. Without loss of generality, we can assume that $A$ is complete, so $A$ is the quotient of a regular local ring. By Lemma \ref{preserve_module}, we can find equations that preserve the $\pd M$, $\dim M$, and $\dim\Ext^{d-r}(M, A)$. By Theorem \ref{reduction}, there is a local ring of characteristic $p$ that satisfies those equations, which produces a counterexample in characteristic $p$ as well, contradicting Theorem \ref{ext_dim_lowpd}.
\end{proof}

Similar arguments can be used to extend Corollary \ref{dim_one_dim} to equicharacteristic, and to reduce the Grade Conjecture from equicharacteristic zero to characteristic~$p$:

\begin{cor}
Let $A$ be a local ring in equicharacteristic of dimension $d$, and suppose $\pd M < \infty$ and $\dim M = 1$. Then $\dim \Ext^{d-1}(M, A) = 1$.
\end{cor}

\begin{thm}
Suppose the Grade Conjecture is true for all local rings of characteristic $p$. Then it is true for all equicharacteristic local rings.
\end{thm}

\textbf{Acknowledgments.} This paper is based in part on my dissertation, and I would like to thank my advisor, Sankar Dutta, for his advice and support throughout my graduate career. I would also like to thank the referee for comments that substantially improved the paper.

\end{document}